\documentclass[a4paper]{article}

\usepackage[english]{babel}
\usepackage[utf8x]{inputenc}
\usepackage[T1]{fontenc}

\usepackage[a4paper,top=3cm,bottom=2cm,left=3cm,right=3cm,marginparwidth=1.75cm]{geometry}

\usepackage{amsmath}
\usepackage{amssymb}
\usepackage{graphicx}
\usepackage[colorinlistoftodos]{todonotes}
\usepackage[colorlinks=true, allcolors=blue]{hyperref}
\usepackage{enumerate}   
\newtheorem{proposition}{Proposition}
\newtheorem{theorem}[proposition]{Theorem}
\newtheorem{lemma}[proposition]{Lemma}
\newtheorem{corollary}[proposition]{Corollary}
\newtheorem{definition}[proposition]{Definition}

\newtheorem{remark}[proposition]{Remark}

\numberwithin{equation}{section}

\numberwithin{proposition}{section}

\usepackage{color}

\author{T. Bloom, L. Bos, N. Levenberg, S. Ma'u and F. Piazzon}

\begin{document}

\newcommand{\C}{{\mathbb{C}}}
\newcommand{\R}{{\mathbb{R}}}
\newcommand{\Z}{{\mathbb{Z}}}
\newcommand{\B}{{\mathbb{B}}}

\title{The Extremal Function for the Complex Ball for Generalized Notions of Degree and Multivariate Polynomial Approximation}
\maketitle
\begin{center} {\bf Dedicated to the memory of Professor Jozef Siciak.} \end{center}
\begin{abstract}
We discuss the Siciak-Zaharjuta extremal function of pluripotential theory for the unit ball in $\C^d$ for spaces of polynomials with the notion of degree determined by a convex body $P.$ We then use it to analyze the approximation properties of such polynomial spaces, and how these may differ depending on the function $f$ to be approximated.
\end{abstract}

\section{Introduction}
The classical Bernstein-Walsh theorem relates the order of approximation of an analytic function in terms of its analyticity inside of level sets of the Siciak-Zaharyuta extremal function. Specifically
 
\begin{theorem} \label{basic} Let $K\subset \C^d$ be compact, nonpluripolar with $V_{K}$ continuous. 
Let $R > 1$, and let $\Omega_R := \{ z : V_{K} (z) < \log R \}$.
Let $f$ be continuous on $K$. Then
$$
  \limsup_{n\to \infty} D_n(f,K)^{1/n}\leq 1/R
$$ 
if and only if $f$ is
the restriction to $K$ of a function holomorphic in $\Omega_R$.  
\end{theorem}

Here for $K\subset \C^d$ compact,
\begin{equation}\label{vkfcn}V_K(z) = \max[0,\sup \{\frac{1}{deg (p)}\log |p(z)|: ||p||_K:=\max_{\zeta \in K}|p(\zeta)| \leq 1\}],\end{equation}
where $p$ is a nonconstant holomorphic polynomial, is the Siciak-Zaharyuta extremal function for $K$ and for a continuous 
complex-valued function $f$ on $K,$
\[ D_n(f,K) := \inf \{ ||f-p_n||_{K}:p_n \in \mathcal P_n \}\]
is the error in best uniform approximation to $f$ on $K$ by polynomials of degree at most $n$. We write ${\mathcal P}_n$ for the space of holomorphic polynomials of degree at most $n$. 
 
Recently Trefethen \cite{T} has argued that polynomial approximation on the hypercube $K=[-1,1]^d\subset \R^d$ by the space of polynomials of what he refers to as of {\it euclidean degree} at most $n$ can be quite advantageous. By this is meant the space of polynomials
\[\{p\in\R[x]\,\, p(x)=\sum_{|\alpha|_2\le n}a_\alpha x^\alpha,\,\,x\in\R^d,\,\,a_\alpha\in\R\}\]
where for the multi-index $\alpha\in\Z_+^d,$ ${|\alpha|_2:=\sqrt{\sum_{i=1}^d\alpha_i^2}}$ is the usual euclidean norm of $\alpha.$

Generalizations of the notion of the degree of a polynomial and the associated extremal functions have been given by Bayraktar \cite{Bay}. Indeed, given a convex body $P\subset  (\R^+)^d=[0,\infty)^d$ we may define a {\it $P-$extremal function} $V_{P,K}$ associated to $K$. Specifically,
we suppose that $P\subset (\R^+)^d$ is a compact convex set in $(\R^+)^d$ with non-empty interior $P^o$. 
We also require that $P\subset (\R^+)^d$ has the property that 
\begin{equation}\label{phyp} \Sigma \subset kP \ \hbox{for some} \ k\in \Z^+\end{equation}
where
\[\Sigma:=\{(x_1,...,x_d)\in \R^d: x_1,...,x_d \geq 0, \ x_1+\cdots x_d\leq 1\}\]
is the standard (unit) simplex.

Associated with $P$, following \cite{Bay}, we consider the finite-dimensional polynomial spaces 
\[Poly(nP):=\{p(z)=\sum_{J\in nP\cap (\Z^+)^d}c_J z^J: c_J \in \C\}\]
for $n=1,2,...$. Here $J=(j_1,...,j_d)$. In the case $P=\Sigma$ we have $Poly(n\Sigma)=\mathcal P_n$, the usual space of holomorphic polynomials of degree at most $n$ in $\C^d$. 

Another class of  examples is given by
$P_q:=\{(x_1,...,x_d)\in (\R^+)^d: (x_1^q+\cdots x_d^q)^{1/q}\leq 1\}$, the (nonnegative) portion of an $l^q$ ball in $(\R^+)^d$, $1\leq q\leq \infty$.

Note that $P_1=\Sigma$ and hence $Poly(nP_1)=\mathcal P_n$ while $P_2$ is that part of the euclidean ball in the positive "octant" and so $Poly(nP_2)$ corresponds to the space of polynomials of "euclidean degree" at most $n$ considered by Trefethen.

Clearly there exists a minimal positive integer $A=A(P)\geq 1$ such that $P\subset A\Sigma$. Thus
\begin{equation}\label{AP} Poly(nP) \subset Poly(An\Sigma)=\mathcal P_{An} \ \hbox{for all} \ n.\end{equation}
We let $d_n:={\rm dim}(Poly(nP))$ and note that by (\ref{AP}), $d_n=O(n^d).$ It follows from convexity of $P$ that 
\[ p_n \in Poly(nP), \ p_m \in Poly(mP) \Rightarrow p_n \cdot p_m\in Poly((n+m)P).\]

Now, recall the {\it indicator function} of a convex body $P$ is
\[\phi_P(x_1,...,x_d):=\sup_{(y_1,...,y_d)\in P}(x_1y_1+\cdots + x_dy_d).\]
For the $P$ we consider, $\phi_P\geq 0$ on $(\R^+)^d$ with $\phi_P(0)=0$.  
Define the logarithmic indicator function 
\[H_P(z):=\sup_{J\in P} \log |z^J|:=\phi_P(\log |z_1|,...,\log |z_d|).\]
Here $|z^J|:=|z_1|^{j_1}\cdots |z_d|^{j_d}$ for $J=(j_1,...,j_d)\in P$ (the components $j_k$ need not be integers). From (\ref{phyp}), we have
\[H_P(z)\geq \frac{1}{k}\max_{j=1,...,d}\log^+ |z_j|=\frac{1}{k}\max_{j=1,...,d}[\max(0,\log |z_j|)].\]
We use $H_P$ to define generalizations of the Lelong classes $L(\C^d)$, the set of all plurisubharmonic (psh) functions $u$ on $\C^d$ with the property that $u(z) - \log |z| = 0(1), \ |z| \to \infty$, and 
\[L^+(\C^d)=\{u\in L(\C^d): u(z)\geq \log^+|z| + C_u\}\]
where $C_u$ is a constant depending on $u$. We remark that, a priori, for a set $E\subset \C^d$, one defines the {\it global extremal function}
\[V_E(z):=\sup \{u(z): u\in L(\C^d), \ u\leq 0 \ \hbox{on} \ E\}.\]
It is a theorem, due to Siciak and to Zaharjuta (cf., Theorem 5.1.7 in \cite{K}), that for $K\subset \C^d$ compact, $V_K$ coincides with the function in (\ref{vkfcn}). Moreover, 
\[V_K^*(z):=\limsup_{\zeta\to z}V_K(\zeta)\in L^+(\C^d)\]
precisely when $K$ is {\it nonpluripolar}; i.e., for $K$ such that $u$ plurisubharmonic on a neighborhood of $K$ with $u=-\infty$ on $K$ implies $u\equiv -\infty$. 

Define
\[L_P=L_P(\C^d):= \{u\in PSH(\C^d): u(z)- H_P(z) =0(1), \ |z| \to \infty \},\] and 
\[L_{P,+}=L_{P,+}(\C^d)=\{u\in L_P(\C^d): u(z)\geq H_P(z) + C_u\}.\]
Then $L_{\Sigma} = L(\C^d)$ and $L_{\Sigma,+} = L^+(\C^d)$. Given $E\subset \C^d$, the {\it $P-$extremal function of $E$} is given by $V^*_{P,E}(z):=\limsup_{\zeta \to z}V_{P,E}(\zeta)$ where
\[V_{P,E}(z):=\sup \{u(z):u\in L_P(\C^d), \ u\leq 0 \ \hbox{on} \ E\}.\]
For $P=\Sigma$, we recover $V_E=V_{\Sigma,E}$. We will restrict to the case where $E=K\subset \C^d$ is compact. In this case, Bayraktar \cite{Bay} proved a Siciak-Zaharjuta type theorem showing that $V_{P,K}$ can be obtained using polynomials. Note that $\frac{1}{n}\log |p_n|\in L_P$ for $p_n\in Poly(nP)$. 

\begin{proposition} \label{turgay2} Let $K\subset \C^d$ be compact and nonpluripolar. Then 
\[V_{P,K} =\lim_{n\to \infty} \frac{1}{n} \log \Phi_n\]
pointwise on $\C^d$ where
\[\Phi_n(z):= \sup \{|p_n(z)|: p_n\in Poly(nP),  \ ||p_n||_K\leq 1\}.\]
If $V_{P,K}$ is continuous, the convergence is locally uniform on $\C^d$.
\end{proposition}

\noindent Note that $V_{P,K}=0$ on the polynomial hull $\hat K$ of $K$. Also, if $V_K$ is continuous, so is $V_{P,K}$ (cf., the discussion after Proposition 2.3 in \cite{BL17}).

The degree of approximation of analytic functions by polynomials in $Poly(nP)$ is given by a generalization of the Bernstein-Walsh Theorem proved in \cite{BL17}. With the notation 
\[
  D_n = D_n(f,K,P) 
  := 
  \inf \{ ||f-p_n||_{K}:p_n \in Poly(nP)\},
\]
\begin{theorem} (\cite{BL17}) \label{BWthm} Let $K$ be compact and assume $V_{P,K}$ is continuous. Let $R > 1$, and let $\Omega_R=\Omega_{R(P,K)} := \{ z : V_{P,K} (z) < \log R \}$. Let $f$ be continuous on $K$. 
Then $f$ is
the restriction to $K$ of a function holomorphic in $\Omega_{R(P,K)}$ if and only if
\[
  \limsup_{n\to \infty} D_n(f,P,K)^{1/n}\leq 1/R.
\]
\end{theorem}

Reference \cite{BL17} also gives a formula for the $P-$extremal function of a product set.  We make the following definition: we call a convex body $P\subset (\R^+)^d$ a {\it lower set} if for each $n=1,2,...$, whenever $(j_1,...,j_d) \in nP\cap (\Z^+)^d$ we have $(k_1,...,k_d) \in nP\cap (\Z^+)^d$ for all $k_l\leq j_l, \ l=1,...,d.$ 

\begin{proposition}(\cite{BL17}) \label{productproperty} Let $P\subset (\R^+)^d$ be a lower set and let $E_1,...,E_d\subset \C$ be compact and nonpolar. Then
\begin{equation}\label{prodprop} V^*_{P,E_1\times \cdots \times E_d}(z_1,...,z_d)=\phi_P(V^*_{E_1}(z_1),...,V^*_{E_d}(z_d)).\end{equation}
\end{proposition}

They use this formula to explain the (sometimes) advantageous approximation properties of polynomial spaces of euclidean degree at most $n$ discovered by Trefethen.

In this work we discuss the case of $K=B_2:=\{z\in\C^d\,:\, \|z\|_2\le1\},$ the complex unit ball in $\C^d,$ as an example of a non-product set. Based on the approach discussed in the next section, we get an explicit formula for $V_{P_\infty,B_2}$ in Proposition \ref{ellinf}. We analyze the approximation properties on $K=B_2$ of polynomial spaces $Poly(nP_q)$ as in Theorem \ref{BWthm} and see how these may differ depending on the function $f$ in section 3. In section 4 we compute the Monge-Amp\`ere measure $(dd^cV^*_{P,K})^2=\mu_{P_{\infty},B_2}$ (Proposition \ref{MAmeas}) and give a probabilistic application following \cite{Bay}. 

The genesis of this work took place at the Dolomites Research Week in Approximation, September 4-8, 2017. 

\section{Computing extremal functions}

To compute extremal functions, in particular $V_{P,B_2}$ for various $P$, we will generalize the approach of Bloom \cite{B97}, for which we will require a generalized version of a theorem of Zeriahi \cite{Ze} (see also \cite[Theorem 3.2]{B97}) that allows one to compute the extremal function by means of orthogonal polynomials. Hence consider $K\subset\C^d$ a compact set and let $\mu$ be a finite Borel measure supported on $K$ satisfying a Bernstein-Markov inequality, i.e., for every $\epsilon>0$ there exists a constant $C(\epsilon)>0$ such for all  holomorphic polynomials $p\in\C[z],$
\begin{equation}\label{BM}
\|p\|_K\le C(\epsilon)(1+\epsilon)^{{\rm deg}(p)}\|_{L^2(\mu)}.
\end{equation}
Here ${\rm deg}(p)$ denotes the usual degree of $p.$ However, associated to the polyhedron $P$ we may define
\begin{definition}
For a holomorphic polynomial $p\in\C[z]$ ($z\in\C^d$), we set
\[{\rm deg}_P(p):=\inf \{n\in\Z^+\,:\, p\in Poly(nP)\}\]
and for $\alpha\in(\Z^+)^d,$
\[|\alpha|_P:=\inf \{t\in\R^+\,:\, z^\alpha\in Poly(tP)\}, \]
i.e., the classical Minkowski norm of the vector $\alpha$ with respect to $P.$\end{definition}

We note that
\begin{equation}\label{degs}
{\rm deg}_P(z^\alpha)-1\le |\alpha|_P\le {\rm deg}_P(z^\alpha). 
\end{equation}

We remark that by our assumption (\ref{phyp}) on $P$ we may equivalently replace the classical degree (${\rm deg}(p)={\rm deg}_\Sigma(p)$) in (\ref{BM}) by ${\rm deg}_P(p).$

To define the orthogonal polynomials we impose an ordering on the multinomial indices $\alpha\in (\Z^+)^d,$ which is consistent with the degree, i.e.,
\[\alpha\le \beta \,\,\implies\,\, {\rm deg}_P(z^\alpha)  \le {\rm deg}_P(z^\beta). \]
We then let
\[\{p_\alpha(z)=p_\alpha(z,\mu)\,:\,\alpha\in(\Z^+)^d\}\]
be the family of orthonormal polynomials obtained by the Gram-Schmidt process with inner-product given by $\mu$ applied to the
monomials $\{z^\alpha\,:\,\alpha\in(\Z^+)^d\}$ so ordered.

\begin{theorem} (Generalized Zeriahi \cite{Ze}) \label{zeriahi} Under the above assumptions
\[V_{P,K}(z)=\limsup_{\alpha}\, \frac{1}{|\alpha|_P} \log\left|p_\alpha(z)\right|\quad\hbox{for}\,\,z\in\C^d\setminus \hat{K}\]
where $\hat{K}$ denotes the polynomial hull of $K.$
\end{theorem}
\noindent {\bf Proof}. The argument is a straightforward generalization of that of Zeriahi. We give the details for the sake of completeness.

First note, that by our assumption that $\mu$ satisfies a Bernstein-Markov inequality (\ref{BM}), 
\[\limsup_{\alpha}\frac{1}{{\rm deg}_P(p_\alpha)}\log|p_\alpha(z)|\le V_{P,K}(z),\,\,z\in \C^d\setminus\widehat{K}.\]
Then, as ${\rm deg}_P(p_\alpha)={\rm deg}_P(z^\alpha)$, from (\ref{degs}) we also have
\[\limsup_{\alpha}\frac{1}{|\alpha|_P}\log|p_\alpha(z)|\le V_{P,K}(z),\,\,z\in \C^d\setminus\widehat{K}.\]

To show the reverse inequality, first recall that by Proposition \ref{turgay2} we have
\begin{equation}\label{V_limit}
V_{P,K}(z)=\lim_{n\to\infty}\left(\sup \frac{1}{n} \log|p(z)|\,:\, p\in Poly(nP)\,\,\hbox{and}\,\, \|p\|_K\le1\right).
\end{equation}
Now, let $q\in Poly(nP)$ be such that $\|q\|_K\le 1.$ We expand $q$ in its orthogonal series with repsect to the basis
$\{p_\alpha\,:\, {\rm deg}_P(p_\alpha)\le n\},$ i.e.,
\[q(z)=\sum_{\alpha\in nP} c_\alpha p_\alpha(z)\]
where 
$$
c_\alpha = \int_K q(z)\overline{p_\alpha(z)}d\mu(z).
$$
Since $\|q\|_K\le1$ we have
\[|c_\alpha|\le \int_K |p_\alpha(z)|d\mu(z)\le \sqrt{\mu(K)}\]
by the Cauchy-Schwarz inequality. Thus
\begin{equation}\label{q_bound}
|q(z)|\le {\rm dim}(Poly(nP))\,\sqrt{\mu(K)}\,\max_{\alpha\in nP}|p_\alpha(z)|=d_n\sqrt{\mu(K)}\,\max_{\alpha\in nP}|p_\alpha(z)|.
\end{equation}
Now fix a $z_0\in \C^d\setminus \widehat{K}$ and let $\alpha_n$ be the largest multiindex in our ordering such that
\[|p_{\alpha_n}(z_0)|=\max_{\alpha\in nP} |p_\alpha(z_0)|.\]
We note that by the fact that the chosen ordering respects ${\rm deg}_P$ we have that $n\le m$ implies that ${\rm deg}_P(z^{\alpha_n})\le {\rm deg_P}(z^{\alpha_m}),$ i.e., the sequence $\{{\rm deg}_P(z^{\alpha_n})\}$ is monotonically increasing. Further,
the sequence of multi-indices $\{\alpha_n\}$ satisfies $\lim_{n\to\infty}{\rm deg}_P(z^{\alpha_n})=+\infty$ for, if not, say 
${\rm deg}_P(z^{\alpha_n})\le M$ for all $n,$ then, by (\ref{q_bound}) for any polynomial $q(z)$ satisfying $\|q\|\le 1,$
we have
\[|q(z_0)|\le d_n\sqrt{\mu(K)}\,\max_{\alpha\in MP} |p_\alpha(z_0)|\]
so that by (\ref{V_limit}) $V_{P,K}(z_0)=0,$ a contradiction.

Thus we also have $\lim_{n\to\infty} |\alpha_n|_P=+\infty$ and, using (\ref{V_limit}) and (\ref{q_bound}), we have
\[V_{P,K}(z_0)\le \limsup_{n\to\infty} \frac{1}{n} \log|p_{\alpha_n}(z_0)|.\]
But, as noted previously, $|\alpha_n|_P\le {\rm deg}_P(z^{\alpha_n})\le n$ and so also
\[\limsup_{\alpha}\, \frac{1}{|\alpha|_P} \log\left|p_\alpha(z)\right|\ge V_{P,K}(z_0).\]
$\square$

For $K=B_2$ normalized monomials will be used in Theorem \ref{zeriahi} in the next section. It is worth noting that, for slightly more general compact sets $K$, the monomials are also Chebyshev polynomials. Specifically, 
consider $K\subset \mathbb{C}^d$ compact. Let $<_l$ be the lexicographic ordering on the multiindices $\alpha \in (\mathbb{Z}^+)^d$ given by $\alpha>_l\beta$ if $|\alpha|>|\beta|$ or
if $|\alpha|=|\beta|$ and $\alpha_i=\beta_i$ for $i=1,...,r$ and $\alpha_{r+1}>\beta_{r+1}$ for some $r$.

For each multiindex $\alpha$ we define a collection  $\mathcal{Q}(\alpha)$ of polynomials as follows. Let \[\mathcal{Q}(\alpha):=\{q(z)=z^{\alpha} +\sum_{\beta <_l\alpha}c_{\beta}z^{\beta}: c_{\alpha}\in \mathbb{C}\}.\]
Let $b_{\alpha}=\inf_{q\in Q(\alpha)}||q||_K$ and let \[B_{\Sigma}=\{\theta\in(\mathbb{R}^+)^d\quad:\quad|\theta|=1\}\]
and we write $B_{\Sigma}^{\circ}=\{\theta\in B_{\Sigma}: \theta_i>0, \ i=1,...,d\}$ for its interior.

The following result is due to Zaharjuta \cite{Za}.
\begin{theorem}
For $\theta\in B_{\Sigma}^{\circ}$ 
we have
\[b(\theta,K):=\lim_{\frac{\alpha}{|\alpha|}\rightarrow\theta}b_{\alpha}^{\frac{1}{|\alpha|}}\] exists and
$\log b(\theta,K)$ is convex on $B_{\Sigma}^{\circ}$.
\end{theorem}

\noindent The number $b(\theta,K)$ is called the {\it directional Chebyshev constant} (with direction $\theta$) for $K.$

For $\mu$ a Bernstein-Markov measure on $K$ (cf. \eqref{BM}) we set 
\[h_{\alpha}:=\inf_{q\in \mathcal{Q}(\alpha)}||q||_{L^2(\mu)}.\]
Then we have

\begin{proposition}
For $\theta\in B_{\Sigma}$
\[\lim_{\frac{\alpha}{|\alpha|}\rightarrow\theta}b_{\alpha}^{\frac{1}{|\alpha|}}= \lim_{\frac{\alpha}{|\alpha|}\rightarrow\theta}h_{\alpha}^{\frac{1}{|\alpha|}}\] 
in the sense that one of the limits exists if and only if the other does and in that case both are equal.
\end{proposition}

We say that a polynomial $q_0$ \textit{realizes} $b_{\alpha},$ i.e., $q_0$ is a Chebyshev polynomial for $K$ of  index $\alpha,$ if
$q_0\in \mathcal{Q}(\alpha)$ and $||q_0||_K=b_{\alpha}$ (and similarly for $h_{\alpha}$). Now assume that $K$ is invariant under the torus action
\[z=(z_1,...,z_d)\rightarrow(e^{it_1}z_1,...,e^{it_d}z_d), \ t_1,...,t_d\in \R\]
and that $\mu$ is also invariant under the torus action. This will be the case in the next section. Then the monomials are mutually orthogonal and any polynomial which realizes $h_{\alpha}$ is a monomial. Moreover, we have
    
\begin{proposition}
Let $K$ be invariant under the torus action. For each multiindex $\alpha$ the monomial $z^{\alpha}$ realizes $b_{\alpha},$ i.e., $z^\alpha$ is a Chebyshev polynomial for $K.$
\end{proposition} 
\noindent{\bf Proof}.
Let $q_0$ be a polynomial which realizes $b_{\alpha}.$ 
Suppose that $\alpha=(\alpha_1,...,\alpha_d)$ and $\alpha_1>0$. Let 
$$q_1(z):=\frac{1}{\alpha_1}\sum_{j=1}^{\alpha_1}q_0(e^{\frac{2\pi ij}{\alpha_1}}z_1,z_2,...,z_d).$$
Then $q_1$ is homogeneous in $z_1$ of degree $\alpha_1$, $q_1\in \mathcal{Q}(\alpha)$, and $$||q_1||_K\leq||q_0||_K$$ so $||q_1||_K=b_{\alpha}.$

Then repeat successively the averaging procedure for each of the remaining variables $z_j$ for which $\alpha_j>0$. We obtain the monomial $z^{\alpha}$ and we see that $||z^{\alpha}||_K=b_{\alpha}.$
$\square$

Note that for $K$ invariant under the torus action and $P$ a convex set, the polynomial spaces $Poly(nP)$ as well as the extremal function $V_{P,K}(z)$ are invariant under the torus action.

\section{The case of $K$ the unit ball in $\C^d$}
Here we take
\[K=B_d:=\{z\in\C^d\,:\, |z|\le1\}\]
where $|z|:=\|z\|_2=(|z_1|^2+\cdots |z_d|^2)^{1/2}$ denotes the euclidean norm of $z\in\C^d$ and $\mu$ denotes Lebesgue measure on $K$. It is well-known that (\ref{BM}) holds in this setting.

It is known (see e.g. \cite{R}) that the monomials $z^\alpha$ are mutually orthogonal and indeed
\[p_\alpha(z)=c_\alpha z^\alpha,\,\, \alpha\in(\Z^+)^d,\]
with
\[c_\alpha^2:=\frac{(|\alpha|+d)!}{\alpha! \pi^d}\]
are the orthonormal polynomials. Here $|\alpha|:=\sum_{j=1}^d\alpha_j$ and $\alpha!:=\prod_{j=1}^d (\alpha_j!).$

Now, for $0\neq z\in\C^d$ 
let
\[I(z):=\{i\,:\, z_i\neq 0\}.\]
By Theorem \ref{zeriahi}, in this case the extremal function is determined by the limsup of the sequence of normalized monomials. However, by compactness, any sequence of normalized multiindices $\alpha(j)/|\alpha(j)|_P$ has a limit point, and hence we first consider such convergent sequences.

\begin{lemma}\label{BloomFunctional}
Suppose that $\{\alpha(j)\in(\Z^+)^d\}$ is an infinite sequence of {\it distinct} multi-indices, ordered as above, such that
 $i\not\in I(z)\,\,\implies\,\,\alpha_i(j)=0,$ and that
\[\lim_{j\to\infty} \frac{\alpha(j)}{|\alpha(j)|_P}=\theta\in(\R^+)^d.\]
Necessarily then $|\theta|_P=1$ and $i\notin I(z)\,\,\implies\,\,\theta_i=0.$

We have
\begin{align*}
\lim_{j\to\infty} \frac{1}{|\alpha(j)|_P}\log|c_{\alpha(j)}z^{\alpha(j)}|&=F_d(\theta; z)\\
&:=\frac{1}{2}\left\{
\sum_{i\in I(z)}\theta_i\log(|z_i|^2)-\sum_{i\in I(z)}\theta_i\log(\theta_i)+\left(\sum_{i\in I(z)}\theta_i\right)
\log\left(\sum_{i\in I(z)}\theta_i\right)\right\}.
\end{align*}
\end{lemma}
\noindent {\bf Proof}. This is a straightforward calculation based on Stirling's formula, $\log(m!)=m\log(m)-m+O(\log(m))$ and the fact that, by construction, $\lim_{j\to\infty} |\alpha(j)|_P=\infty.$
$\square$

\begin{proposition} \label{prop22} For $z\in \C^d\setminus K$
\[V_{P,K}(z)=\max_{\theta\in (\R^+)^d, \, |\theta|_P=1} F_d(\theta; z).\]
\end{proposition}
\noindent {\bf Proof}. First note that restricted to any hyperplane of the form $\{w\in\C^d\,:\,w_i=0\}$
the unit ball, the extremal function and, at least for points $z$ such that $i\notin I(z)$, the functional $F_d(\theta;z)$ all reduce to the same corresponding lower dimensional problem. Hence we may, without loss of generality, assume that $z_i\neq0,$ $1\le i\le d,$ i.e., $I(z)=\{1,2,\ldots,d\},$ and
\[F_d(\theta;z)=\frac{1}{2}\left\{
\sum_{i=1}^d\theta_i\log(|z_i|^2)-\sum_{i=1}^d\theta_i\log(\theta_i)+\left(\sum_{i=1}^d\theta_i\right)
\log\left(\sum_{i=1}^d\theta_i\right)\right\}.\]
The proof is now straightforward as $K=B_d$ is polynomially convex, by Theorem \ref{zeriahi} we have
\[V_{P,K}(z)=\limsup_{\alpha}\, \frac{1}{|\alpha|_P} \log\left|c_\alpha z^\alpha\right|\quad\hbox{for}\,\,z\in\C^d\setminus K.\]
Further, every convergent subsequence of $\{\alpha/|\alpha_P|\}$ has its limit in $\{\theta\in (\R^+)^d\,:\, |\theta|_P=1\}$ and
every such $\theta$ is the limit of such a subsequence. Combined with Lemma \ref{BloomFunctional} the result follows.
$\square$

\medskip For the sake of completeness we will now verify the known formula for the extremal function in the  case of the classical degree, i.e., when $P=\Sigma,$ the standard unit simplex.

\begin{proposition}\label{ClassicalExtremal} For $P=\Sigma,$
\begin{enumerate}[(i)]
\item \[\max_{\theta\in (\R^+)^d\,:\, \sum_{i=1}^d \theta_i=1}F_d(\theta:z)=\log(|z|),\,\, \forall \,z\in\C^d,\,z\neq0,\]
\item \[V_{K}(z):=V_{\Sigma,K}(z)=\log(|z|),\,\,z\notin K.\]
\end{enumerate}
\end{proposition}
\noindent {\bf Proof}. Formula (ii) follows immediately from (i). To show (i) we proceed by induction on the dimension. When $d=1,$ $\theta_1=1$ and trivially
\[F_1(\theta;z)=\frac{1}{2}\left\{ 1\times \log(|z|^2)-0\right\}=\log(|z|).\]

We suppose then that the result holds for up to dimension $d-1$ and must prove that it also holds for dimension $d.$ Again, we may assume without loss that $z_i\neq0,$ $1\le i\le d.$ We maximize
\begin{align*}
F_d(\theta;z)&=\frac{1}{2}\left\{
\sum_{i=1}^d\theta_i\log(|z_i|^2)-\sum_{i=1}^d\theta_i\log(\theta_i)+\left(\sum_{i=1}^d\theta_i\right)
\log\left(\sum_{i=1}^d\theta_i\right)\right\}\\
&=\frac{1}{2}\left\{
\sum_{i=1}^d\theta_i\log(|z_i|^2)-\sum_{i=1}^d\theta_i\log(\theta_i)\right\}
\end{align*}
over the set $\{\theta\in (\R^+)^d\,:\, \sum_{i=1}^d \theta_i=1\}.$ 

Consider first the interior ($\theta_i>0,\,\,\forall i$) critical point(s) given by Lagrange multipliers as the solution of
\[\log(|z_i|^2)-(1+\log(\theta_i))=\lambda,\,\,1\le i\le d,\]
or equivalently,
\begin{align*}
\log(|z_i|^2)-\log(\theta_i)&=\log(|z_d|^2)-\log(\theta_d),\,\,1\le i\le d\\
\iff \log(|z_i|^2/\theta_i)&=\log(|z_d|^2/\theta_d), \,\,1\le i\le d\\
\iff |z_i|^2/\theta_i&=|z_d|^2/\theta_d, \,\,1\le i\le d\\
\iff \theta_i&=\theta_d\frac{|z_i|^2}{|z_d|^2}, \,\,1\le i\le d.
\end{align*}
Taking the sum of both sides we see that
\[1=\sum_{i=1}^d \theta_i = \theta_d\frac{1}{|z_d|^2}\sum_{i=1}^d |z_i|^2=\theta_d\frac{1}{|z_d|^2}|z|^2\]
so that 
\[\theta_d=\frac{|z_d|^2}{|z|^2}\]
and
\[\theta_i=\theta_d\frac{|z_i|^2}{|z_d|^2}=\frac{|z_i|^2}{|z|^2}, \,\,1\le i\le d. \]
Substituting these values of the $\theta_i$ into the expression for $F_d$ we obtain the critical value of
\begin{equation}\label{crit1}
F_d(\theta;z)=\frac{1}{2}\left\{\sum_{i=1}^d \frac{|z_i|^2}{|z|^2}\log(|z_i|^2)-\sum_{i=1}^d \frac{|z_i|^2}{|z|^2}\log\left(\frac{|z_i|^2}
{|z|^2}\right)\right\}
=\log(|z|),
\end{equation}
after simplfication.

The other competitors for the maximum are on the boundary of our constraint set $\{\theta\in (\R^+)^d\,:\, \sum_{i=1}^d \theta_i=1\},$ i.e., when one or more of the $\theta_i$ are equal to zero. But in this case we reduce to a lower dimensional version of the same problem, and by our induction assumption the maximum of $F_d(\theta;z)$ is then
\[\frac{1}{2} \log\left(\sum_{i\,:\,\theta_i\neq0}|z_i|^2\right)\]
which is less than the value at the interior critical point. Hence the maximum is indeed $\log(|z|)$ and we are done. $\square$

\medskip
We next collect some basic facts about the function $F_d.$
\begin{proposition}\label{FBasics}
Assume again that $z_i\neq0,$ $1\le i\le d.$ Then
\begin{enumerate}[(i)]
\item $F_d(\theta;z)$ is homogeneous of order one in $\theta$ so that
\[F_d(\theta;z)=\sum_{i=1}^d \theta_i \frac{\partial}{\partial \theta_i}F_d(\theta;z);\]
\item $\nabla F_d(\theta;z):=(\frac{\partial F_d}{\partial \theta_1},...,\frac{\partial F_d}{\partial \theta_d})\neq 0\in\R^d$ for $z\notin K;$
\item At any interior point, $\theta_i>0,$ $1\le i\le d,$ the Hessian of $F_d(\theta;z)$ is non-positive definite;
\item If $|z|<1$ then $F_d(\theta;z)<0;$ 
\item If $|z|=1$ then $\displaystyle{\max_{\theta\in (\R^+)^d,\,|\theta|_P=1}F_d(\theta;z)=0};$
\item If $|z|>1$ then $\displaystyle{\max_{\theta\in (\R^+)^d,\,|\theta|_P=1}F_d(\theta;z)>0}.$
\end{enumerate}
\end{proposition}
\noindent {\bf Proof}. Item (i) is completely elementary and so we leave out the details. For (ii) we calculate
\[2\frac{\partial}{\partial \theta_i}F_d(\theta;z)=\log(|z_i|^2)-\log(\theta_i)+\log\left(\sum_{i=1}^d\theta_i\right) \]
so that  $\nabla  F_d(\theta;z)= 0$ iff 
\begin{align*}
\log(|z_i|^2)&=\log\left(\theta_i/\left(\sum_{j=1}^d\theta_j\right)\right),\,\,1\le i\le d\\
\iff |z_i|^2&=\theta_i/\left(\sum_{j=1}^d \theta_j\right), \,\,1\le i\le d
\end{align*}
and hence, taking the sum, we must have
\[\sum_{i=1}^d |z_i|^2=\left(\sum_{i=1}^d\theta_i\right)/\left(\sum_{j=1}^d\theta_j\right)=1.\]
To see (iii), we easily calculate
\[2\frac{\partial^2}{\partial \theta_i\partial\theta_j}F_d(\theta;z)=\begin{cases}\frac{1}{S} -\frac{1}{\theta_i}&{\rm if}\,\,j=i\\
\frac{1}{S}&{\rm if}\,\,j\neq i \end{cases}
\]
where $\displaystyle{S:=\sum_{i=1}^d\theta_i}.$ Hence (twice the) Hessian, $H_F,$ say, is
\[2H_F=\frac{1}{S}u u^t-D\]
where 
\[u:=
\left[\begin{array}{c}1 \cr 1\cr \cdot \cr \cdot\cr 1\end{array} \right]\in \R^d\,\,{\rm and}\,\, D:= 
\left[\begin{array}{ccccc} 1/\theta_1&0&\cdot&\cdot&0 \cr 0&1/\theta_2&0&\cdot&0\cr
0&\cdot&\cdot&\cdot&0\cr
0&\cdot&\cdot&0&1/\theta_d\end{array}\right]\in\R^{d\times d}\]
which we recognize as a rank one perturbation of the negative definite diagonal matrix $-D.$ More specifically, it is easy to
verify that
\[H_Fv=0\in \R^d\,\,\hbox{where}\,\,v:=\left[\begin{array}{c}\theta_1\cr\theta_2\cr\cdot\cr\cdot\cr\theta_d\end{array}\right]\]
and that for $0\neq w\in\R^d$ with $w^tu=0,$
\[w^t(2H_F)w=-w^tDw<0\]
so that $H_F$ is singular and negative definite on a $(d-1)$-dimensional subspace of $\R^d.$

Properties (iv), (v) and (vi) can be easily verified using the homogeneity. Indeed, for any $\theta\in (\R^+)^d$ with $|\theta|_P=1,$ there is a $\theta'\in (\R^+)^d$ with $\sum_{i=1}^d\theta'_i=1$ and $t_\theta>0$ such that $\theta=t_\theta\theta'$ and hence
\[F_d(\theta;z)=F_d(t_\theta\theta';z)=t_\theta F_d(\theta';z)\]
and the result follows from the classical case, Proposition \ref{ClassicalExtremal}. $\square$

For brevity's sake let
\[B_P:=\{\theta\in(\R^+)^d\,:\,|\theta|_P=1\}\]
denote the constraint set. 

\begin{lemma}\label{bdry} 
Suppose that $B_P$ is a smooth manifold near its boundary. Then, if  $z_i\neq0,$ $1\le i\le d,$ the maximum of $F_d(\theta;z)$ over $B_P$
is never attained at a boundary point, i.e., where one or more of the $\theta_i=0.$
\end{lemma}
\noindent {\bf Proof}. We just note that as
\[2\frac{\partial}{\partial \theta_i}F_d(\theta;z)=\log(|z_i|^2)-\log(\theta_i)+\log\left(\sum_{i=1}^d\theta_i\right) \]
it follows that 
\[\lim_{\theta_i\to0^+}\frac{\partial}{\partial \theta_i}F_d(\theta;z)=+\infty \]
while the partials otherwise are finite. Hence a sufficiently small positive perturbation of $\theta_i=0$ will result in an increase in the value of $F_d(\theta;z).$ $\square$

\begin{remark} Note this means, e.g., that
for such $P$ we can never have $V_{P,B_2}(z_1,z_2)=V_{P,B_1}(z_1)$ for a point $(z_1,z_2)\in\C^2$ with $z_2\neq0.$ 
\end{remark}

\begin{lemma}\label{uniqueness} Suppose that $B_P$ is strictly convex (i.e. $|(x+y)/2|_P<(|x|_P+|y|_P)/2,$ $x\neq y$). Then  if  $z_i\neq0,$ $1\le i\le d,$ and $z\not\in K,$ the maximum of $F_d(\theta;z)$ over $B_P$ is uniquely attained.
\end{lemma}
\noindent {\bf Proof}. Suppose for the sake of a contradiction that the maximum is attained at two distinct points 
$\theta',\theta''\in B_P.$ By Lemma \ref{bdry} both $\theta',\theta''$ are in the interior of $B_P.$ Now, by
Proposition \ref{FBasics}, (iii), $F_d(\theta;z)$ is a concave function so that 
\[F_d((\theta'+\theta'')/2 ;z)\ge(F_d(\theta';z)+F_d(\theta'';z))/2=\max_{\theta\in B_P} F_d(\theta;z). \]
Note that, as $|(\theta'+\theta'')/2|_P<(|\theta'|_P+|\theta''|_P)/2=1,$ 
\[t:=2/|\theta'+\theta''|_P>1.\]
Then, $\theta:=t(\theta'+\theta'')/2\in B_P$ and
\[F_d(\theta;z)=F_d(t(\theta'+\theta'')/2;z)=tF_d((\theta'+\theta'')/2;z)\geq t\max_{\theta\in B_P} F_d(\theta;z)>
\max_{\theta\in B_P} F_d(\theta;z),\]
a contradiction.
$\square$

\medskip In case the $P-$norm is a smooth function, the maximum can be characterized by Lagrange multipliers. For simplicity's sake
let $g(\theta):=|\theta|_P$ so that $B_P=\{\theta\in(\R^+)^d\,:\,g(\theta)=1\}.$ Then the Lagrange multiplier equations are
\[\frac{\partial}{\partial \theta_i}F_d(\theta;z)=\lambda \frac{\partial}{\partial \theta_i}g(\theta)\,\,1\le i\le d.\]
Taking the sum of both sides we obtain by homogeneity
\[ F_d(\theta;z)=\sum_{i=1}^d\theta_i\frac{\partial}{\partial \theta_i}F_d(\theta;z)
=\lambda \sum_{i=1}^d \theta_i \frac{\partial}{\partial \theta_i}g(\theta). \]
But as $g(\theta)$ is a norm, it also is homogeneous of order one, and so by the Euler identity, the sum on the righthand side reduces to 
$g(\theta)=1$ for $\theta\in B_P.$ In other words, the Lagrange multiplier
\begin{equation}\label{multiplier}
\lambda =  F_d(\theta;z).
\end{equation}

\subsection{The case of $P$ the unit $\ell_q$ ball, $1\le q\le \infty$}
\begin{figure}
\centering
\includegraphics[width=1.0\textwidth]{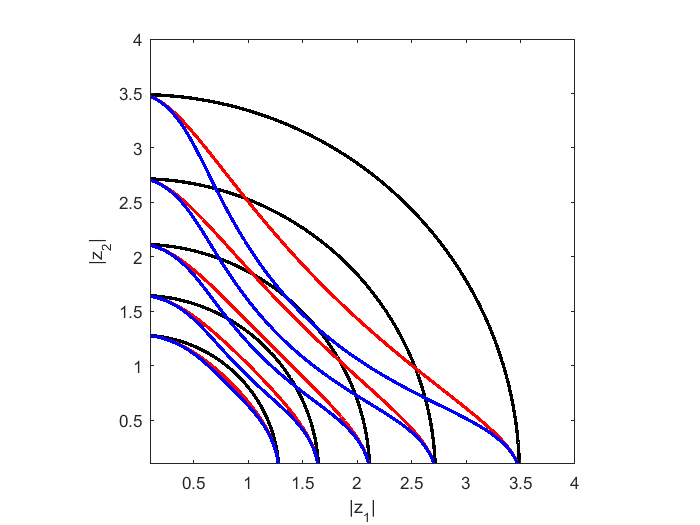}
\caption{\label{fig:fig1}Extremal Function Contour Plots for Levels $.25, .5, .75, 1, 1.25.$ The black curves
correspond to $q=1,$ the red ones to $q=2$ and the blue ones to $q=4.$}
\end{figure}

For $1\leq q <\infty,$
\[g(\theta)=|\theta|_P=\left\{\sum_{i=1}^d\theta_i^q\right\}^{1/q}\]
is smooth and hence the associated extremal function may be found by solving the Lagrange multipliers equations
\begin{align*}
\frac{\partial}{\partial \theta_i}F_d(\theta;z)&=\lambda \frac{\partial}{\partial \theta_i}g(\theta)\\
&=F_d(\theta;z)\frac{\theta_i^{q-1}}{g(\theta)^{q-1}}\\
&=F_d(\theta;z)\theta_i^{q-1},\,\,1\le i\le d
\end{align*}
using \eqref{multiplier}. These $d$ equations are actually dependent as the sum of both sides multiplied by $\theta_i$ reduces to the tautology $F_d(\theta;z)=F_d(\theta;z).$ Hence we solve the system
\begin{align*}
\frac{\partial}{\partial \theta_i}F_d(\theta;z)&=F_d(\theta;z)\theta_i^{q-1},\,\,1\le i\le (d-1)\\
g(\theta)&=1.
\end{align*}
For $z_i\neq 0,$ $1\le i\le d,$ a unique solution is guaranteed by Lemma \ref{uniqueness}. In the case of $d=2$ this is particularly easy to find numerically and in Figure \ref{fig:fig1} we show several contours for $q=1,2,4.$ One notices immediately that on the diagonal $|z_1|=|z_2|$ the extremal functions are notably different, whereas they have the same values on the complex lines $z_1=0$ and $z_2=0$ (as here they all reduce to the same univariate extremal function; cf., the proof of Proposition \ref{prop22}).

This has some interesting consequences for the approximation of functions from $Poly(nP).$  Let $P_q$ denote the unit $\ell_q$ ball intersect $(\R^+)^d.$

In the following three examples $K$ denotes the euclidean unit ball $B_2$ in $\C^2$.

\medskip\noindent {\bf Example 1}. Consider
\[f_1(z_1,z_2):=\frac{1}{1-z_1/2}+\frac{1}{1-z_2/2}.\]
As the monomials form an orthogonal basis, best $L_2$ approximations are equivlent to Taylor expansions. In this case we have
\[f_1(z_1,z_2)=\sum_{k=0}^n\left(\frac{z_1^k}{2^k}+\frac{z_2^k}{2^k}\right)+\left\{\frac{(z_1/2)^{n+1}}{1-z_1/2}+
\frac{(z_2/2)^{n+1}}{1-z_2/2}\right\}\]
for $|z_1|, |z_2|<2$, in particular, on $K$. 
But note that for {\it any} $q\ge1,$ the degrees of $z_1^k$ and $z_2^k$ are both $k.$ In particular, the best $L_2$ approximation for $f_1$ on $K$ of degree $n,$ for {\it any} $q\ge1$ is
\[p_n(z_1,z_2):=\sum_{k=0}^n\left(\frac{z_1^k}{2^k}+\frac{z_2^k}{2^k}\right).\]
In other words, for $L_2$ approximation, there is no advantage in a higher value of $q$ despite the fact that the spaces $Poly(nP_q)$ are of increasing dimension in $q.$
$\square$

\medskip Less obvious examples may be analyzed by means of the extremal function.
Indeed, by Theorem \ref{BWthm} the order of uniform approximation to a holomorphic function $f(z)$ by $Poly(nP)$ is given (essentially)
by
\[D_n(f,K,P)=O(R^{-n})\]
where
\[\log(R):=\inf_{z\in S(f)} V_{P,K}(z)\]
and $S(f)\subset \C^d$ is the singular set of $f.$

\medskip\noindent {\bf Example 2}. Consider the bivariate Runge type function
\[f_2(z_1,z_2):=\frac{1}{a^2+z_1^2+z_2^2}\quad a>1.\]
Its singular set is given by
\[S(f_2):=\{z\in \C^2\,:\, z_1^2+z_2^2=-a^2\}.\]
\begin{lemma} We have
\[\min_{z\in S(f_2)}V_{P_q,K}(z)=\log(a),\,q\ge1,\]
attained at (among other points) $z_1=ia, z_2=0.$
\end{lemma}
\noindent {\bf Proof}.
Consider first the classical case $q=1$ when $P_q=\Sigma.$ By Proposition \ref{ClassicalExtremal} the extremal function 
is $\log(|z|).$ Hence it suffices to  show that
\[\min_{z_1^2+z_2^2=-a^2}|z_1|^2+|z_2|^2=a^2.\]
To see this we calculate
\begin{align*}
\min_{z_1^2+z_2^2=-a^2}|z_1|^2+|z_2|^2&=\min_{z_1\in\C}|z_1|^2+|a^2+z_1^2|^2\\
&=\min_{r\ge0\,\,\theta\in[0,2\pi]} r^2+\{a^4+2a^2r^2\cos(2\theta)+r^4\}^{1/2}\quad \hbox{(writing } z_1=r\exp(i\theta) )\\
&= \min_{r\ge0}\,\, r^2+\{a^4-2a^2r^2+r^4\}^{1/2}\quad (\hbox{for }\theta=\pi/2)\\
&=\min_{r\ge0} \,\,r^2+\{(a^2-r^2)^2\}^{1/2}\\
&=\min_{r\ge0} \,\,r^2+|a^2-r^2|\\
&=\min_{r\ge0}\,\,\begin{cases} a^2&\hbox{if }r\le a\cr 
2r^2-a^2&\hbox{if } r\ge a\end{cases}\\
&=a^2.
\end{align*}
A particular minimum point is given by $r=a,\theta=\pi/2,$ i.e., $z_1=ia$ for which $z_2^2=-a^2-z_1^2=0.$

For any other value of $\infty\ge q>1,$ we note that $Poly(nP_q)\supset Poly(nP_1)$ and hence the approximation error
\[D_n(f_2,K,P_q)\le D_n(f_2,K,P_1)\] and so comparing the orders of error decay we must have
\[\min_{z\in S(f_2)}V_{P_q,K}(z)\ge \min_{z\in S(f_2)}V_{P_1,K}(z)=\log(a).\]
On the other hand $V_{P_q,K}(ia,0)=V_{P_1,K}(ia,0)=\log(a)$ and so also
\[\min_{z\in S(f_2)}V_{P_q,K}(z)\le\log(a).\]
$\square$

\begin{figure}
\centering
\includegraphics[width=1.0\textwidth]{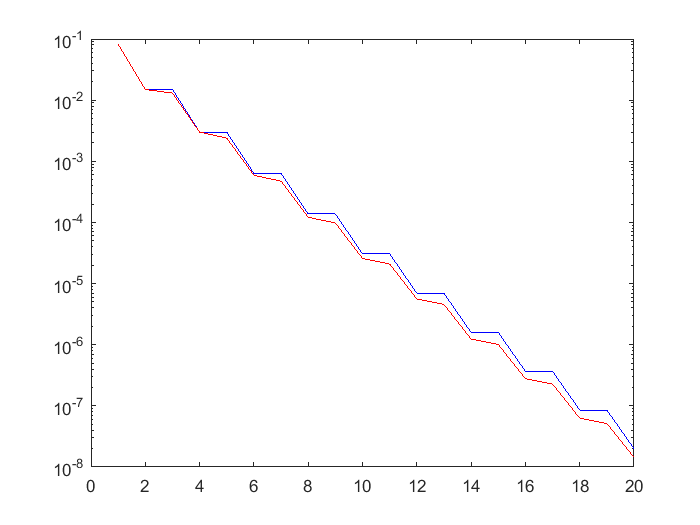}
\caption{\label{fig:fig2}Errors in best $L_2$ approximation to $f_2$; blue curve for $q=1$ and red (lower) curve for $q=4$}
\end{figure}

\medskip In other words the rate of decay of the uniform approximation errors to $f_2$ are also the same for all choices of $q\ge1;$ there is no approximation value added despite the fact that the dimensions of the spaces $Poly(nP_q)$ are increasing in $q$. This behavior is illustrated numerically in Figure \ref{fig:fig2} where we show the $L_2$ best approximation error for $f_2$ with $a=2$ as a function of $n$ for $q=1$ and $q=4.$ $\square$

\medskip\noindent{\bf Example 3}. There is no gain in approximating $f_1$ or $f_2$ by the spaces $Poly(nP_q),$ $q>1,$ precisely because there is a singular point on the coordinate hyperplane $z_2=0$ where the extremal functions all reduce to the same univariate extremal function for all $q\ge1.$ We now give an example of a function whose singular set does not approach the coordinate hyperplanes and for which the approximation order of $Poly(nP_q)$ is strictly increasing in $q.$ Specifically, let
\[f_3(z_1,z_2)=\frac{1}{1-z_1z_2}.\]
The best $L_2$ approximation is again easy to calculate by means of a Taylor series, which in this case is just a geometric series:
\[f_3(z_1,z_2)=\frac{1}{1-z_1z_2}=\sum_{k=0}^m z_1^kz_2^k + \frac{(z_1z_2)^{m+1}}{1-z_1z_2}\]
for $|z_1z_2|<1$, in particular, on $K$. 
The {\it uniform} norm of the error on $K$ is easily bounded by
\begin{equation}\label{UnifErr}
\max_{|z|\le1} \left| \frac{(z_1z_2)^{m+1}}{1-z_1z_2}\right|\le \frac{2^{-(m+1)}}{1-1/2}=2^{-m}.
\end{equation}
If we take $m=n/2$ (ignoring round-offs) then we approximate $f_3$ by a polynomial
\[p_n(z_1,z_2):=\sum_{0\le k\le n/2} z_1^kz_2^k\]
of classical degree $n.$ Its uniform error is then $O(2^{-n/2})$ implying that
\[\min_{z\in S(f_3)} V_{P_1,K}(z)\ge \log(\sqrt{2}).\]
On the other hand, for $z_0=(1,1)\in S(f_3),$
\begin{align*}V_{P_1,K}(z_0)&\ge F_2(\theta;z_0)\,\,{\rm with}\,\,\theta=(1/2,1/2)\\
&=\frac{1}{2}\left\{
\frac{1}{2}\log(1)+\frac{1}{2}\log(1)-(\frac{1}{2}\log(\frac{1}{2})+\frac{1}{2}\log(\frac{1}{2}))+(\frac{1}{2}+\frac{1}{2})
\log(1)\right\}\\
&=\log(\sqrt{2})
\end{align*}
so that also
\[\min_{z\in S(f_3)} V_{P_1,K}(z)\le \log(\sqrt{2})\]
and we may conclude that
\[\min_{z\in S(f_3)} V_{P_1,K}(z)= \log(\sqrt{2})\]
and that the rate of decay of the uniform error, $O(2^{-n/2}),$ is optimal for $q=1.$

For values of $q>1,$ note that $|(k,k)|_q\le n$ iff $k\le n/2^{1/q}.$ Hence
\[p_n(z):=\sum_{0\le k\le n/2^{1/q}}z_1^kz_2^k \in Poly(nP_q)\]
with uniform error on $K$ of $O(2^{-n/2^{1/q}}=O((2^{(2^{-1/q})})^{-n})$ by \eqref{UnifErr}. Again, this implies that
\[\min_{z\in S(f_3)} V_{P_q,K}(z)\ge \log(2^{(2^{-1/q})}).\]
On the other hand, again  for $z_0=(1,1)\in S(f_3),$ 
\begin{align*}V_{P_q,K}(z_0)&\ge F_2(\theta;z_0)\,\,{\rm with}\,\,\theta=(2^{-1/q},2^{-1/q})\\
&=\log(2^{(2^{-1/q})})
\end{align*}
and we may conclude, that in general,
\[\min_{z\in S(f_3)} V_{P_q,K}(z)=\log(2^{(2^{-1/q})}) \]
and the optimal rate of decay of the uniform error is $O((2^{(2^{-1/q})})^{-n}).$ For example, as $q\to\infty$ this rate approaches $O(2^{-n}),$ considerably better than the $O(2^{-n/2})$ for the classical case. Note also that this advantage persists even when the difference in the dimensions of the various polynomial spaces is taken into account. Indeed, for the classical total degree the dimension of the bivariate polynomials of degree at most $n$ is $N:=(n+2)(n+1)/2$ so that the decay of the error in terms of the dimension is $O(2^{-n/2})=O(2^{-\sqrt{N/2}}).$ On the other hand, for the tensor-product case, $q=\infty,$ the dimension is $N=(n+1)^2$ so that the error decays like $O(2^{-n})=O(2^{-\sqrt{N}}).$ $\square$

\medskip

We do not believe that there is a closed formula for the extremal function for $1<q<\infty.$ However for $q=\infty$ we may show that
\begin{proposition}\label{ellinf} Suppose that $d=2.$ Then for $|z|\ge1,$
\[V_{P_\infty,B_2}(z)=\begin{cases}
\frac{1}{2}\left\{\log(|z_2|^2)-\log(1-|z_1|^2)\right\}&\hbox{if  } |z_1|^2\le 1/2\,\,{\rm and}\,\,|z_2|^2\ge1/2 \cr
\frac{1}{2}\left\{\log(|z_1|^2)-\log(1-|z_2|^2)\right\}&\hbox{if  } |z_1|^2\ge 1/2\,\,{\rm and}\,\,|z_2|^2\le1/2 \cr
\log(|z_1|)+\log(|z_2|)+\log(2)&\hbox{if  } |z_1|^2\ge 1/2\,\,{\rm and}\,\,|z_2|^2\ge1/2 \end{cases}.
\]
\end{proposition}
\noindent {\bf Proof}. If $z_1=0$  the first case of the formula reduces to $\log(|z_2|),$ i.e., the univariate extremal function in $z_2,$ as is correct. Similarly, if $z_2=0$ the second case of the formula reduces to $\log(|z_1|),$ i.e., the univariate extremal function in $z_1,$ as is correct. Hence we suppose that $z_1,z_2\neq0$ and we maximize $F_d(\theta;z)$ over the constraint 
\[B_{P_\infty}=\{(\theta_1,\theta_2)\,:\,0\le\theta_1,\theta_2\le1\}.\]
Lemma \ref{bdry} informs us that the maximum cannot be attained at a boundary point of $B_{P_\infty},$ i.e., when
either $\theta_1,\theta_2=0.$

Consider first the upper edge of the constraint $\theta_2=1,$ $0\le \theta_1\le 1.$ The boundary value at $\theta_1=\theta_2=1$
\begin{equation}\label{corner}
F_d((1,1));z)=\frac{1}{2}\left\{\log(|z_1|^2)+\log(|z_2|^2)-0+2\log(2)\right\}=\log(|z_1|)+\log(|z_2|)+\log(2)
\end{equation}
is a candidate for the maximum (while, as mentioned above $\theta_1=0, \ \theta_2=1$ is not). Competitors are given by critical points along this edge. Hence we calculate
\begin{align*}&\frac{\partial}{\partial \theta_1}F_d(\theta_1,1;z)=\frac{1}{2}\left\{\log(|z_1|^2)-\log(\theta_1)+\log(\theta_1+1)\right\}=0\\
\iff &\log(|z_1|^2)=\log(\theta_1)-\log(\theta_1+1)=\log(\theta_1/(\theta_1+1))\\
\iff & |z_1|^2=\theta_1/(\theta_1+1)\\
\iff & \theta_1= |z_1|^2/(1-|z_1|^2) \quad (|z_1|\neq1). 
\end{align*}
Now it is easy to check that $\theta_1= |z_1|^2/(1-|z_1|^2)\in[0,1]$ iff $|z_1|^2\le 1/2,$ i.e., we have a competitor critical point in this case and otherwise we do not.
If indeed, $|z_1|^2\le1/2$ then we calculate
\begin{align*}
F_d((|z_1|^2/(1-|z_1|^2),1);z)&=\frac{1}{2}\left\{
\frac{|z_1|^2}{1-|z_1|^2}\log(|z_1|^2)+\log(|z_2|^2)- \frac{|z_1|^2}{1-|z_1|^2}\log\left(\frac{|z_1|^2}{1-|z_1|^2}\right)\right.\\
&\quad+\left.\left(\frac{|z_1|^2}{1-|z_1|^2}+1\right)\log\left(\frac{|z_1|^2}{1-|z_1|^2}+1\right)\right\}\\
&=\frac{1}{2}\left\{ \log(|z_2|^2)-\log(1-|z_1|^2)\right\}
\end{align*}
after some simplification. 

Now, we claim that this critical value, in the case that $|z_1|^2\le1/2,$ is greater than the corner value \eqref{corner}.
Indeed, 
\begin{align*}
&\frac{1}{2}\left\{ \log(|z_2|^2)-\log(1-|z_1|^2)\right\}\ge \log(|z_1|)+\log(|z_2|)+\log(2)\\
\iff & \frac{1}{2}\left\{ \log(|z_2|^2)-\log(1-|z_1|^2)\right\}\ge \frac{1}{2}\left\{ \log(|z_1|^2)+\log(|z_2|^2)+\log(4)\right\}\\
\iff & -\log(1-|z_1|^2)\ge \log(|z_1|^2)+\log(4)\\
\iff & \log(4|z_1|^2(1-|z_1|^2))\le 0  \\
\iff & 4|z_1|^2(1-|z_1|^2)\le 1\\
\end{align*}
which clearly holds. In summary, we have shown that
\[\max_{0\le\theta_1\le 1,\,\theta_2=1}F_d(\theta;z)=
\begin{cases}
\frac{1}{2}\left\{ \log(|z_2|^2)-\log(1-|z_1|^2)\right\}&\hbox{if  } |z_1|^2\le 1/2 \cr 
\log(|z_1|)+\log(|z_2|)+\log(2)&\hbox{if  } |z_1|^2\ge 1/2
\end{cases}.\]

We immediately obtain the maximum value on the right edge by symmetry, i.e.,
\[\max_{0\le\theta_2\le 1,\,\theta_1=1}F_d(\theta;z)=
\begin{cases}
\frac{1}{2}\left\{ \log(|z_1|^2)-\log(1-|z_2|^2)\right\}&\hbox{if  } |z_2|^2\le 1/2 \cr 
\log(|z_1|)+\log(|z_2|)+\log(2)&\hbox{if  } |z_2|^2\ge 1/2
\end{cases}\]
and the result follows.
$\square$


\section{Computing the extremal measure} Returning to our general setting of a compact, nonpluripolar compact set $K\subset \C^d$ and a convex body $P\subset (\R^+)^d$, recall that $d_n$ is the dimension of $Poly(nP)$. We write
$$Poly(nP)= \hbox{span} \{e_1,...,e_{d_n}\}$$ 
where $\{e_j(z):=z^{\alpha(j)}\}_{j=1,...,d_n}$ are the standard basis monomials. For 
points $\zeta_1,...,\zeta_{d_n}\in \C^d$, let
$$VDM(\zeta_1,...,\zeta_{d_n}):=\det [e_i(\zeta_j)]_{i,j=1,...,d_n}  $$
$$= \det
\left[
\begin{array}{ccccc}
 e_1(\zeta_1) &e_1(\zeta_2) &\ldots  &e_1(\zeta_{d_n})\\
  \vdots  & \vdots & \ddots  & \vdots \\
e_{d_n}(\zeta_1) &e_{d_n}(\zeta_2) &\ldots  &e_{d_n}(\zeta_{d_n})
\end{array}
\right]$$
and for a compact subset $K\subset \C^d$ let
$$V_n =V_n(K):=\max_{\zeta_1,...,\zeta_{d_n}\in K}|VDM(\zeta_1,...,\zeta_{d_n})|.$$
Points $z_1^{(n)},...,z_{d_n}^{(n)}\in K$ achieving the maximum are called Fekete points of order $n$ for $K,P$. It was shown in \cite{BBL} that the limit
$$ \delta(K):=\delta(K,P):= \lim_{n\to \infty}V_{n}^{1/l_n} $$ exists where 
$$l_n:=\sum_{j=1}^{d_n} {\rm deg}(e_j)= \sum_{j=1}^{d_n} |\alpha(j)|$$ is the sum of the degrees of a set of these basis monomials for $ Poly(nP)$. The quantity $\delta(K)$ is called the {\it $P-$transfinite diameter} of $K$. One of the key results in \cite{BBL} was the following:

\begin{theorem} \label{asympwtdfek} Let $K\subset \C^d$ be compact and nonpluripolar. For each $n$, take points $z_1^{(n)},z_2^{(n)},\cdots,z_{d_n}^{(n)}$ in $K$ for which 
$$
 \lim_{n\to \infty}|VDM(z_1^{(n)},\cdots,z_{d_n}^{(n)})|^{1 \over  l_n}=\delta(K)
$$
({\it asymptotically} $P-$Fekete arrays) and let $\mu_n:= \frac{1}{d_n}\sum_{j=1}^{d_n} \delta_{z_j^{(n)}}$. Then
$$
\mu_n \to \frac{1}{d! Vol(P)}(dd^cV^*_{P,K})^d \ \hbox{weak}-*.
$$
\end{theorem}

\noindent Here $vol(P)$ denotes the $\R^d-$Lebesgue measure of $P$.

This shows the significance in being able to find the ``target'' measure $\mu_{P,K}:=(dd^cV^*_{P,K})^d$. It is important to observe that $\mu_{P,K}$ has support in $K$. In this section, we begin with calculations of $\mu_{P,K}$ for certain $P$ and the unit $d-$torus in $\C^d$ and then we use the calculations in the previous section to compute $\mu_{P,K}$ for certain $P$ and the unit ball in $\C^2$. We first recall two results (cf., \cite{Bay} or \cite{BBL}).

For $P\subset (\R^+)^d$ a convex body and $K=T^d$, the unit $d-$torus in $\C^d$, we have  
$$ V_{P,T^d}(z)= H_P(z)=\max_{J\in P} \log(|z^J|)\in L_P^+.$$
If $P=\Sigma=P_1$, then $V_{P,T^d}(z)=\max_{j=1,...,d} \log^+(|z_j|)$. Let $\omega:=dd^c  \max_{j=1,...,d} \log^+(|z_j|)$. We normalize so that $\int_{\C^d} \omega^d=1$. Then for any $u\in L_P^+$ we have
\begin{equation}\label{norm}\int_{\C^d} (dd^cu)^d =\int_{\C^d} (dd^c H_P)^d = d! Vol(P). \end{equation}
where $Vol(P)$ denotes the euclidean volume of $P\subset (\R^+)^d$. In particular, $\mu_{P,K}(K)=d! Vol(P)$.

For simplicity, we take $d=2$; i.e., we work in $\C^2$ and start with $T=\{(z_1,z_2):|z_1|=|z_2|=1\}$. We know that 
$$V_{P,T}(z_1,z_2) = H_P(\log^+|z_1|,\log^+|z_2|).$$ 
Then $V_{P_1,T}(z_1,z_2) = \max[\log^+|z_1|,\log^+|z_2|]$ and $\mu_{P_1,T}$ is normalized Haar measure on $T$. Note that $\mu_{P,T}(T)=2Vol(P)=1$. At the other extreme, for $P_{\infty}=[0,1]\times [0,1]$,
$$V_{P_{\infty},T}(z_1,z_2) = \log^+|z_1|+\log^+|z_2|=\max[0,\log|z_1|,\log|z_2|, \log|z_1|+\log|z_2|].$$ 
We see that near the face $|z_1|=1, \ |z_2|<1$, $V_{P,T}(z_1,z_2) =\log^+|z_1|$ which is maximal there ($(dd^c\log^+|z_1|)^2=0$); ditto for the face $|z_2|=1, \ |z_1|<1$. Thus, as we knew, $\mu_{P_{\infty},T}$ is supported in $T$ but the total mass is $2\cdot Vol([0,1]\times [0,1])=2$. Indeed, for any $1\leq q\leq \infty$, we have 
$$V_{P_q,T}(z_1,z_2) = [(\log^+|z_1|)^{q'}+(\log^+|z_2|)^{q'}]^{1/q'}$$
where $1/q +1/q'=1$. By invariance under $(z_1,z_2)\to (e^{i\theta_1}z_1,e^{i\theta_2}z_2)$, $\mu_{P_q,T}$ is a multiple of normalized Haar measure on $T$; precisely, $\mu_{P_q,T}(T)=2Vol(P_q)$. 

We now turn to the case of the closed Euclidean ball and $P_{\infty}$. We have shown that
\begin{equation}\label{Vformula}
V_{P_\infty,B_2}(z)=\begin{cases}
\frac{1}{2}\left\{\log(|z_2|^2)-\log(1-|z_1|^2)\right\}&\hbox{if  } |z_1|^2\le 1/2\,\,{\rm and}\,\,|z_2|^2\ge1/2 \cr
\frac{1}{2}\left\{\log(|z_1|^2)-\log(1-|z_2|^2)\right\}&\hbox{if  } |z_1|^2\ge 1/2\,\,{\rm and}\,\,|z_2|^2\le1/2 \cr
\log|z_1|+\log|z_2|+\log(2)&\hbox{if  } |z_1|^2\ge 1/2\,\,{\rm and}\,\,|z_2|^2\ge1/2 \end{cases}.
\end{equation}

\begin{proposition} \label{MAmeas} For $K=B_2$ and $P=P_{\infty}$, 
the measure $(dd^cV^*_{P,K})^2=\mu_{P_{\infty},B_2}$ is Haar measure on the torus $\{|z_1|=1/\sqrt 2, \ |z_2|=1/\sqrt 2\}$ with total mass $2$.  
\end{proposition}

\noindent {\bf Proof}.
The function $\frac{1}{2}[\log(|z_1|^2) - \log (1- |z_2|^2)]$ is pluriharmonic in a neighborhood $U_a$ of any point  $a\in\partial B_2\cap \{|z_2|>1/\sqrt{2}>|z_1|\}$, and negative on $U_a \cap B_2$.  Utilizing Proposition 3.8.1 of \cite{K}, we conclude that on $U_a$, $V_{P_{\infty},B_2}= \max\{\frac{1}{2}[\log(|z_1|^2) - \log (1- |z_2|^2)],0\}$ is maximal. Similarly, $V_{P_{\infty},B_2}$ is maximal in a neighborhood of any point of $\partial B_2\cap\{|z_1|>1/\sqrt{2}>|z_2|\}$.  Thus there is no Monge-Amp\`ere mass on these portions of $\partial B_2$; we only have mass on the torus in $\partial B_2$ where $|z_1|,|z_2|= 1/\sqrt 2$.  Invariance of (\ref{Vformula}) under $(z_1,z_2)\mapsto(e^{i\theta_1}z_1,e^{i\theta_2}z_2)$ yields that $(dd^cV^*_{P,K})^2$ is a multiple of Haar measure on this torus. The total mass is $2$ by (\ref{norm}) since the volume of $P$ is $1$.
$\square$

Thus for $P=P_1$, $\mu_{P_1,B_2}$ is supported on the entire topological boundary of $B_2$ while  $\mu_{P_{\infty},B_2}$ is supported on a torus. 

As an interesting application, let $\{p_{\alpha}(z)=c_{\alpha}z^{\alpha}\}$ be the orthonormal polynomials for Lebesgue measure on $K=B_2$ as in section 2. Following \cite{Bay}, we can consider, given $P$, random $Poly(nP)$ polynomials of the form $P_n(z)=\sum_{\alpha\in nP}a_{\alpha}p_{\alpha}(z)$ where the coefficients $a_{\alpha}$ are independent, identically distributed complex-valued random variables. For simplicity, we assume that they are complex Gaussian random variables with distribution
$$\phi(t)dm(t)=\frac{1}{\pi}e^{-|t|^2}dm(t)$$
where $dm$ denotes Lebesgue measure on $\C$. We really want to consider random polynomial mappings $F_n(z)=(P_n(z),Q_n(z))$. Thus we get a probability measure $Prob_n$ on $\mathcal F_n$, the random polynomial mappings with $P_n,Q_n\in Poly(nP)$. We can identify $\mathcal F_n$ with $\C^{d_n}\times \C^{d_n}$. Given $F_n\in \mathcal F_n$, let 
$$\tilde Z_{F_n} := (dd^c \frac{1}{n}\log |F_n|)^2=(dd^c [\frac{1}{2n}\log (|P_n|^2+|Q_n|^2)])^2.$$ For generic $F_n$, $\tilde Z_{F_n}$ is, up to a constant, the normalized zero measure on the (finite) zero set $\{P_n=Q_n=0\}$. The expectation ${\bf E}(\tilde Z_{F_n})$ is a measure on $\C^2$ defined, for $\psi\in C_c(\C^2)$, as
$$\bigl({\bf E}(\tilde Z_{F_n}), \psi\bigr)_{\C^2}:=\int_{\C^{d_n}\times \C^{d_n}}(\tilde Z_{F_n},\psi)_{\C^2}\ dProb_n=\frac{1}{\pi^{2d_n}} \int_{\C^{d_n}\times \C^{d_n}} (\tilde Z_{F_n},\psi)_{\C^2}e^{-\sum_{\alpha\in nP}|a_\alpha|^2}\prod_{\alpha\in nP}dm(a_{\alpha})$$ where $(\tilde Z_{F_n},\psi)_{\C^2}$ denotes the action of the measure $\tilde Z_{F_n}$ on $\psi$. In this setting, Bayraktar proved that 
$$\lim_{n\to \infty} {\bf E}(\tilde Z_{F_n})= (dd^c V_{P,K})^2.$$
as measures. Forming the product probability space of sequences of random polynomial mappings
$$\mathcal P:=\otimes_{n=1}^{\infty} (\mathcal F_n,Prob_n)= \otimes_{n=1}^{\infty} (\C^{d_n}\times \C^{d_n},Prob_n),$$
almost surely (a.s.) in $\mathcal P$ we have 
$$\frac{1}{n}\log |F_n|=\frac{1}{2n}\log (|P_n|^2+|Q_n|^2)\to V_{K,P}(z)$$
pointwise in $\C^2$ and in $L^1_{loc}(\C^2)$. Moreover, a.s. in $\mathcal P$ we have 
$$(dd^c \frac{1}{n}\log |F_n|)^2=(dd^c [\frac{1}{2n}\log (|P_n|^2+|Q_n|^2)])^2\to (dd^c V_{P,K})^2.$$
as measures. 

\begin{corollary} With $K=B_2$, for 
\begin{enumerate}
\item $P=P_1 =\Sigma$, ${\bf E}(\tilde Z_{F_n})\to \mu_{P_1,B_2}$, normalized surface area measure on $\partial B_2$; while for
\item $P= P_{\infty}$, ${\bf E}(\tilde Z_{F_n})\to \mu_{P_{\infty}, B_2}$, a multiple of Haar measure on the torus $\{|z_1|=|z_2|=1/\sqrt 2\}$
\end{enumerate}
with analogous statements for the a.s. results.

\end{corollary}


\bibliographystyle{alpha}
\bibliography{sample}

\begin{thebibliography}{BBL17}

\bibitem[Bay17]{Bay}
T.~Bayraktar.
\newblock Zero distribution of random sparse polynomials.
\newblock {\em Mich. Math. J.}, 66(2): 389--419, 2017.

\bibitem[BBL17]{BBL}
T.~Bayraktar, T.~Bloom, and N.~Levenberg.
\newblock Pluripotential theory and convex bodies.
\newblock {\em Math. Sbornik}, DOI:10.1070/SM8893, 2017.

\bibitem[BL17]{BL17}
L.~Bos and N.~Levenberg.
\newblock Bernstein-walsh theory associated to convex bodies and applications
  to multivariate approximation theory.
\newblock {\em CMFT}, to appear: 342--351, 2017.

\bibitem[Blo97]{B97}
T.~Bloom.
\newblock Orthogonal polynomials in $\mathbb{C}^n$.
\newblock {\em Ind. U. Math. Jour.}, 46(2): 427--452, 1997.

\bibitem[Kli93]{K}
M.~Klimek.
\newblock {\em Pluripotential Theory}.
\newblock Oxford U. Press, 1993.

\bibitem[Rud08]{R}
W.~Rudin.
\newblock {\em Function Theory in the Unit Ball of $\mathbb{C}^n$}.
\newblock Springer, Classics of Mathematics, 2008.

\bibitem[Tre17]{T}
N.~Trefethen.
\newblock Multivariate polynomial approximation in the hypercube.
\newblock {\em Proc. Amer. Math. Soc.}, 145(11): 4837--4844, 2017.

\bibitem[Zah75]{Za}
V.~P. Zaharjuta.
\newblock Transfinite diameter, tchebyshev constants, and capacity for compacta
  in $\mathbb{C}^n$.
\newblock {\em Math. USSR Sbornik}, 25: 350--64, 1975.

\bibitem[Zer85]{Ze}
A.~Zeriahi.
\newblock Capacit\'e, constante de chebyshev et polyn\^omes orthogonaux
  associ\'es a un compact de $\mathbb{C}^n$.
\newblock {\em Bull. Sci. Math. (2)}, 109: 325--335, 1985.

\end{thebibliography}

\end{document}